\documentclass[12pt]{amsart}
\usepackage{amsthm,amsmath,amssymb}
\usepackage{addfont}
\usepackage{mathrsfs}
\usepackage[T1]{fontenc}
\numberwithin{equation}{section}

\def\cc{{\mathcal C}}

\def\cq{{\mathcal Q}}

\def\cas{{\mathcal S}}
\def\ct{{\mathcal T}}

\def\bc{{\mathbb C}}

\def\a{\alpha}

\def\g{\gamma}

\theoremstyle{plain}
\newtheorem{lemma}{Lemma}[section]
\newtheorem{proposition}[lemma]{Proposition}
\newtheorem{theorem}[lemma]{Theorem}
\newtheorem{corollary}[lemma]{Corollary}
\theoremstyle{definition}

\newtheorem{definition}[lemma]{Definition}

\begin{document}

\title[Extremal Schur Multipliers   ]{\textsc{  Extremal Schur Multipliers }}

\author[E.~Christensen]{Erik Christensen}
\address{\hskip-\parindent
Erik Christensen, Mathematics Institute, University of Copenhagen, Copenhagen, Denmark.}
\email{echris@math.ku.dk}
\date{\today}
\subjclass[2020]{ Primary: 15A23, 52A21, 81P47. Secondary: 15A60, 46L07, 47B01.}
\keywords{ Schur product, Hadamard product, Schur multiplier, extremal point, convexity, Grothendieck's inequality, completely bounded map, rank of a matrix} 

\begin{center}
{\bf Dedicated to the memory of Bent Fuglede}
\end{center}

\begin{abstract}   The Schur product of two  complex   $ m \times n $ matrices is their entry wise product. We show that an  extremal element $X$  in the convex set of $m \times n $ complex matrices of Schur multiplier norm at most 1 satisfies the inequality rank$(X) \leq \sqrt{m +n}.$ For positive matrices in $M_n(\bc),$ with unit diagonal, we give a characterization of the extremal elements, and show that  such a matrix satisfies rank$(X) \leq \sqrt{n}.$ 
 \end{abstract}

\maketitle

\section{Introduction}
In  Schur's article \cite{Sc} he studies various properties of the entry wise product inside the complex square  matrices indexed by an infinite set. He establishes the result that this  product of 2 positive $n \times n$ matrices is a positive matrix. The other product of two $n \times n$ matrices, which today is named the {\em matrix product, } is by Schur considered to be an analogy to the convolution product, which  we know from harmonic analysis. In the theory of continuous periodic functions of period $2\pi$  we have  the usual point wise product of functions and the convolution product of continuous functions on the unit circle. The latter product  may also be described as the entry wise product of the coefficients in the Fourier series of the functions. In \cite{Ha}, Hadamard studied the product on the Laurent series given as the  entry wise product of their coefficients, and this is the reason why many researchers use the words {\em Hadamard product } for the matrix product, which  we prefer to name after Schur. In the article \cite{C1} we show that the Schur product in $M_n(\bc)$  may be obtained in a way, which is a direct extension  of  Hadamard's product to the setting of functions on the abelian cyclic group  $C_n$ of order $n$ with values in the complex diagonal $n \times n$ matrices. 

The Schur product also plays a role in the study of quantum channels \cite{AS}, and in that respect we do hope that our results may be of use in the theory of quantum information. Our own motivation for studying the Schur product comes from our interest in the subject, which by  Mathematical Reviews is named {\em noncommutative geometry  à la Connes} with the index 58B34. An introduction to that subject may be found in Connes' book \cite{Co}. In Connes' world the action of a differential operator may be described as a Schur multiplication by a complex valued unbounded self-adjoint matrix  indexed by  the integers. 

\smallskip
Around 1980  we  worked on  various aspects of  Grothendieck's work from the fundamental article  \cite{Gr}, so we have been familiar with several versions of inequalities named {\em Grothendieck's inequality, } without  noticing that Grothendieck's inequality for bilinear forms on the space of continuous functions on a compact space has a dual version, which may be expressed in terms of Schur multipliers. The survey article \cite{Pi} by Pisier on the impact of Grothendieck's work describes this dual result in the Theorem  3.2. The {\em Schur multiplier version} of Grothendieck's inequality takes the following form. 

\begin{theorem}
There exists a positive constant $K^\bc_G$ such that for any natural number $n$ and any complex  $n \times n $ matrix  $X$ of 
Schur multiplier norm 1, there exists a natural number $k$ and for each $i,$ with $1 \leq i \leq k$ there exist a positive real $t_i$ 
and vectors $L_i, R_i$ in $\bc^n$ such that \begin{align*} &\|L_i\|_\infty \leq 1, \, \|R_i\|_\infty \leq 1, \, \sum_i t_i \leq K_G^\bc \\ & \forall \, l,m  \in \{1, \dots, n\}  \, : X_{(l,m)} = \sum_i t_i \overline{L_i(l)} R_i(m).\end{align*} 
\end{theorem}

This shows that the matrix $X$ of the theorem is a positive linear combination, bounded by $K^\bc_G,$  of rank 1 matrices of Schur multiplier norm at most 1. It is well known that an extremal Schur multiplier $Y$ in $M_n(\bc)$  of Schur multiplier norm 1 and rank 1 has  the form $Y_{(l,m)} = \overline{L(l)}R(m)$ for a pair of vectors  $L, $ $R$ in $\bc^n,$ such that $|L(l)| = |R(m)| = 1. $  Hence we find, that all the Schur multipliers of norm at most 1 are contained in $K_G^\bc$ times the convex hull of the extremal Schur multipliers of norm 1 and rank 1.   Since we know that $K^\bc_G > 1,$ the convexity theory tells us that the unit ball of the Schur multipliers must have some  extremal  elements of rank more than 1, but very little is, to the best of our knowledge, known about those extremal Schur multipliers of rank more than 1.   
  
In the proofs to come we will work with matrices in $M_{(m,n)}(\bc)$ and linear operators from $\bc^n$ to $\bc^m,$ and we will use the same type of notation for matrices and for operators, but we have tried  only to  use one of the concepts in a single proof.  
  
   We will not go into a more detailed exposition of the applications and  historical  aspects of the Schur product here, but mention that Horn in \cite{Ho} has given a very nice and complete exposition on  the Hadamard product. We will use the book \cite{HJ} by Horn and Johnson as a basic reference for matrix theory and the books \cite{KR} by Kadison and Ringrose as a basic reference for functional analysis and operator theory.

\section{Some properties of the Schur product}
We  will denote the Schur product of two complex $m \times n$ matrices $X$  and $Y$ by $X\circ Y$ and it is defined by $(X\circ Y)_{(i,j)} := X_{(i,j)}Y_{(i,j)}.$  The Schur multiplier $S_X$ acts on  $M_{(m,n)}(\bc)$ by  $S_X(Y):= X \circ Y.$ 
The Schur multiplier norm of $X$ is the norm of $S_X$ as a linear operator on $M_{(m,n)}(\bc)$ when this space is equipped with the operator norm. As far as we know there is no way to directly compute the norm $\|S_X\|$ from the entries in $X,$ except in the case when $X$ is a positive matrix.
 Then it follows from Schur's result on positive matrices and some elementary C*-algebra theory that $$ X \in M_n(\bc)^+ \implies \|S_X\| =  \max\{ X_{(i,i)}: \, 1 \leq i \leq n\,\}.$$
 
 On the other hand it follows from Proposition 3.3 of \cite{Pi},  that the Schur multiplier norm may be described as a factorization norm, and we have given a detailed  exposition of that in \cite{C2}, where we use the concept named {\em column norm of a matrix.} 
 
\begin{definition} Let $X$ be a complex $m \times n $ matrix then for $1 \leq j \leq n$ $X_j$ denotes the $ m \times 1$ column matrix consisting of the $j'$th column in $X$ and  the column norm of $X$ is denoted $\|X\|_c$ and it is given by the formula $\|X\|_c : =  \max\{\|X_j\|_2\, : \, 1 \leq j \leq n\,\}.$ \end{definition}

The following theorem is copied from \cite{C2},  Theorem 1.3 item (iii), and it is, as mentioned above, basically contained in Grothendieck's work  \cite{Gr}, \cite{Pi} Proposition 3.3, but  our presentation focuses on the fact that the objects involved are matrices in $M_{(m,n)}(\bc),$ and our way to the result follows  a  long path, where we use the theory of completely bounded maps on C*-algebras as described by Paulsen in \cite{Pa}, by Smith in \cite{Sm}, by  Livshits in \cite{Li} and by Walter in \cite{Wa}.

\begin{theorem} \label{ScMu}
\begin{itemize}
Let $X$ be a complex $m \times n$ matrix. \item[(i)] Let $L$ be a complex $k \times m $ matrix and $R $ a complex $k \times n$ matrix. If  $X = L^*R$ then $\|S_X\| \leq \|L\|_c\|R\|_c.$ 
\item[(ii)] Let $r$ be the rank of $X,$ there exist matrices $L$ in $M_{(r,m)}(\bc)$ and $R $ in $M_{(r,n)}(\bc)$  such that $\|L\|_c \|R\|_c = \|S_X\|$ and $X = L^*R.$
\end{itemize} 
\end{theorem} 
 
 \smallskip
It should be mentioned that the Schur multiplier norm has the flavour of an $\ell^\infty$ norm in the way that the norm - in many cases -  only represents a  property  of a small part of the matrix. For instance the factorization mentioned in item (ii) above is almost never unique, as the following simple example demonstrates:$$ X:= \begin{pmatrix}
1   \ 0 \\ 0 \ \frac{1}{4}         
\end{pmatrix} = \begin{pmatrix}
1   \ 0 \\ 0 \ \frac{1}{4}         
\end{pmatrix}\begin{pmatrix}
1   \ 0 \\ 0 \ 1         
\end{pmatrix} = \begin{pmatrix}
1   \ 0 \\ 0 \ \frac{1}{2}         
\end{pmatrix}\begin{pmatrix}
1   \ 0 \\ 0 \ \frac{1}{2}         
\end{pmatrix}.
$$

\section{On a full set of vectors} 
We do not know if the concept, we name {\em a full  set of vectors, } exists in the literature already. If so we apologize for our ignorance, and  then we  would like to know the right word for this property. 

\begin{definition}  Let $n, k$ be natural numbers,  $\{\xi_1, \dots, \xi_n\}$ 
be vectors in $\bc^k$ and $P$ the orthogonal projection from $\bc^k$ onto the linear span of the given vectors. 
Let $\omega_{\xi_1}, \dots, \omega_{\xi_n}$ denote the positive  vector functionals on $M_k(\bc)$ given by $\omega_{\xi_j}(Y):= \langle Y\xi_j, \xi_j \rangle.$
We say that the set of vectors  $\{\xi_1, \dots , \xi_n\}$ is full if the null matrix i $M_k(\bc)$  is the only matrix $X$ in $M_k(\bc)$ which satisfies the equations \begin{equation}
X = PXP \text{ and } \forall j \in \{1, \dots, n \}: \omega_{\xi_j}(X) = 0.
\end{equation} 
\end{definition}
 
 Elementary linear algebra  gives a proof of the following proposition/observation  which plays an important role in this article.

 \begin{proposition}  \label{SqDim} 
Let $\{\xi_1, \dots, \xi_n\}$ be a full set of vectors in $\bc^k.$ Then $ \mathrm{dim}\big(\mathrm{span} (\{\xi_1,  \dots, \xi_n\})\big) \leq \sqrt{n}.$
 \end{proposition}
 
 \begin{proof} Let $r$ denote the dimension of the  linear span $\cas$  of the vectors $\{\xi_1, \dots, \xi_n\}, $ and $P$ the orthogonal projection in $M_k(\bc) $ with range $\cas.$ Then the dimension of the space $PM_k(\bc)P$ is $r^2. $ Since the set of vectors is full we get that $n \geq r^2$ and then $r \leq \sqrt{n}.$ 
 \end{proof} 
 
 We have yet another simple observation on the concept {\em fullness.} 
 
 \begin{proposition} \label{TransFull} 
 Let $\{\xi_1, \dots , \xi_n\}$ be a set of vectors $\bc^k$ with linear span $\cas,$ and $T : \cas \to \bc^l $ an injective linear map. The set $\{\xi_1, \dots ,  \xi_n\}$  of vectors in $\bc^k$ is full if and only if the set $\{T\xi_1, \dots ,T\xi_n\}$ of vectors in $\bc^l$ is full. 
 \end{proposition}  
 
 \begin{proof}
 Let $\ct$ denote the linear span of the vectors $\{T\xi_1, \dots ,T\xi_n\},$  then $T$ has an inverse injective linear map of $\ct$ onto $\cas,$ so we will only have to prove, that if the set $\{\xi_1, \dots ,\xi_n\}$ 
is full then the set $\{T\xi_1, \dots ,T\xi_n\}$ is full. Let $Q$ denote the orthogonal projection from $\bc^l $ onto $\ct$ and  $\tilde T$ the linear operator from $\bc^k$  to $\bc^l$ which is defined by the equations $$  \tilde T \g :=  \begin{cases} 0 \quad \text{ if } P\g = 0\\ T\g \, \, \text{ if } P\g = \g. \end{cases} $$ 
Suppose $B$ is a linear operator on $\bc^l$ such that $QBQ = B,$ and for each $j$ we have   $$0 = \langle B\tilde T\xi_j, \tilde T \xi_j\rangle = \langle \tilde T^* B \tilde T \xi_j, \xi_j \rangle. $$ 
Then $\tilde T^* B \tilde T = 0$ and the injectivity of $T$ plus the equation $QBQ = B $ implies that $B = 0,$ so the proposition follows.   
\end{proof}   

\section{Extremal positive Schur multipliers of norm 1} 
We know from above that the positive Schur multipliers in $M_n(\bc) $  of norm 1 are of the form $S_X$ with $X$ positive and $$\max\{X_{(1,1)}, \dots, X_{(n,n)}\} = 1.$$ Let $r$ denote the rank of $X,$ $F$ denote the square root of $X$ and let $W$ denote the matrix of a partial  isometry of the range space of $F$ onto $\bc^r.$ Then for $L$ in $M_{(r,n)}(\bc) $ defined as $L:= WF$ we have $X = L^*L,$ and we have obtained a factorization of $X$ as described in item (ii) of Theorem \ref{ScMu}. This shows that in the case of positive matrices, an optimal Schur factorization from Theorem \ref{ScMu} item (ii)  may be constructed directly, contrary to the general case. 

Before we start to present our result on the extremal positive Schur multipliers, we have to stress that we use the word extremal positive Schur multiplier for those matrices that are the extreme points of the convex set consisting of positive matrices of Schur multiplier norm at most 1. An extremal element here may not be extremal in the unit ball of all Schur multipliers. An example to demonstrate this fact may be seen in the matrix $X$ defined as $X:= \begin{pmatrix} 1 & \frac{1}{2} \\ \frac{1}{2}
 & \frac{1}{4}
\end{pmatrix}. $ This matrix is positive  of rank 1 and any non trivial convex combination of two non proportional positive matrices in $M_2(\bc)$ will have rank 2, so it is an extremal point in the positive matrices of Schur multiplier norm 1. 
On the other hand $X$ is the midpoint between two rank 1 matrices,   $$ X:= \begin{pmatrix} 1 & \frac{1}{2} \\ \frac{1}{2}
 & \frac{1}{4}
\end{pmatrix}  = \frac{1}{2}\bigg(\begin{pmatrix}
1 \\ \frac{1}{2}
\end{pmatrix}\begin{pmatrix}
1 & \frac{3}{4}
\end{pmatrix} +  \begin{pmatrix}
1 \\ \frac{1}{2}
\end{pmatrix}\begin{pmatrix}
1 & \frac{1}{4}
\end{pmatrix}\bigg),$$ 
and by item (i) in Theorem \ref{ScMu}  we see that $X$ is a convex combination of 2 Schur multipliers of norm at most 1.  

\begin{theorem} \label{ExtPos}
Let $X$ be a positive $n \times n$ matrix with square root $F$ and of Schur multiplier norm 1. If the  matrix $X$ is an extremal point  in the set of positive Schur multipliers of norm at most 1,  then  the set of columns from $F$ of norm 1 is a full set of vectors in $\bc^n,$ and the dimension of the linear span of the columns from $F$ of norm 1 is at most $\sqrt{n}.$ 
\end{theorem}

\begin{proof}
Let $J $ denote the set of indices $j$ such that the column $F_j$ has length 1. We will give a proof based on contra position, so suppose  that the set $\cc$  of unit vectors in $\bc^n$ given as $\cc:= \{F_j: j \in J\}$  is not full. Let $P$ denote the matrix of the orthogonal projection from $\bc^n $ onto span$(\cc),$ then there exists  a non zero self-adjoint matrix $B$ in $M_n(\bc)$ such that $B= PBP$ and for all $j $ in $J$ we have $\omega_{F_j} (B) = 0.$ 
Without loss of generality we may replace $B$ by a positive multiple such that we may assume that $B$ is still non trivial and for all $i $ in $\{1, \dots , n \} $ which are not in $J$ we have $\|B\| \leq 1 - \|F_i\|_2^2 \leq 1. $ 
Let I denote the unit in $M_n(\bc),$ then since $\|B\| \leq 1$  we have that both $(I+B) $ and $(I-B)$ are positive matrices in $M_n(\bc)$ so $Y$ and $Z$ defined as $Y:= F(I+B)F$ and $Z:= F(I-B)F$ are positive matrices which satisfy $X = \frac{1}{2}Y + \frac{1}{2}Z.$ 
We will then have to show that the diagonals in $Y$ and $Z$ are of norm at most 1, and that $Y \neq X$ in order to show that $X$ is not extremal. Let us look at the diagonal elements $Y_{(j,j)}.$ If $j$ is in $ J$ then 
$$Y_{(j,j)} = \omega_{F_j} (I) + \omega_{F_j}(B) = 1.$$ If $i$ is not in $J$ then  
$$Y_{(i,i)} = \omega_{F_i} (I) + \omega_{F_i}(B) \leq \|F_i\|_2^2 + (1 - \|F_i\|_2^2) \|F_i\|_2^2 < 1.$$ Since the range space of $F$ contains all the vectors in the set $\cc$ and $0 \neq B = PBP $ we have $FBF \neq 0$ and $X \neq Y,$ so $X$  is not an extremal positive Schur multiplier of norm 1.  

The last statement in the theorem follows directly from Proposition \ref{SqDim}, and the theorem follows. 

\end{proof}

The positive matrices in $M_n(\bc)$ which have the property that their diagonal equals the unit $I$ do play a special role in the applications of Schur multipliers since they posses  the property that for such a matrix $X$ we have that $S_X$ is a completely  positive unit preserving map on  $M_n(\bc),$ and this property is appreciated  by the people studying quantum channels \cite{AS}. We let $\cq_n $ denote the set of  positive $n \times n $ matrices with diagonal $I.$ Contrary to the set of positive Schur multipliers of norm 1, the set $\cq_n$ is a face in the unit ball of the Schur multipliers, as the following proposition shows.  

\begin{proposition} \label{Qface} Let $0 <\a < 1,$   $X$ be in $\cq_n,$ and  $Y,$ $Z$ in $M_n(\bc)$ of Schur multiplier norm at most 1 such that $X = (1- \a)Y + \a Z$ then both $Y$ and $Z$ are in $\cq_n.$  \end{proposition}  

\begin{proof} By Theorem \ref{ScMu} item (ii)we get $Y = A^*B$ and $Z = C^*D$ such that the column norm of each of the matrices $\{A, B, C, D\}$ is at most 1. Since the number 1 is extremal in the unit disk, we have for each index $j$ that  $\langle B_j, A_j \rangle = \langle D_j, C_j \rangle = 1$ so $ A_j = B_j ,$  $C_j = D_j$ and all have length 1. 
We may conclude that $Y=A^*A,$ and $ Z  = C^*C$ with unit diagonals, so they  both belong to $\cq_n,$ and the proposition follows. 
\end{proof}

We can extend Theorem \ref{ExtPos} to the following theorem. 

\begin{theorem} \label{ExtQ} 
A matrix $X$ in $\cq_n$  is extremal in $\cq_n$  if and only if the set of columns $\{X_1, \dots , X_n\}$ is a full set of vectors in $\bc^n.$ An extremal element in $\cq_n$ has rank at most $\sqrt{n}.$ 
\end{theorem}

\begin{proof}
Let $F$ denote the square root of $X,$ then the support  projection of $F$ equals that of $X, $ so the linear map induced by $F$ will act injectively on the linear span of the columns in $F.$ Since $ X= F^2,$ the linear map induced by $F$ will map the columns in $F$ onto the columns in $X$  and we get by Proposition \ref{TransFull}, that the set of columns in $X$ is a full set if and only if the set of columns in $F$ is a full set. We will then show that $X$ is extremal if and only if the columns in $F$ is a full set of vectors. By Theorem \ref{ExtPos} it follows, that if $X$ is extremal then the set of columns in $F$ must be full. Then  let us assume that $X$ is not extremal in $\cq_n$ and  that $X$ is a convex combination, $X = \frac{1}{2}Y  + \frac{1}{2}Z$ of two  matrices from $\cq_n$ with $ Y \neq X.$ As above we let $P$ denote the range projection of $F,$ 
which is also the range projection of $X.$ The equation $X = \frac{1}{2}Y  + \frac{1}{2}Z$ implies that $0 \leq Y \leq 2 X$ so there exists a matrix $C$ in $M_n(\bc)$ such that $\|C\| \leq \sqrt{ 2 }, $ $PCP = C $ and $F \neq \sqrt{Y } = CF.$  Since the theorem on polar decomposition contains  a uniqueness result, we get that  $ C^*C \neq P,$ and we can then define a non zero self-adjoint matrix $B$ by $B:= C^*C- P.$ Then  we have $F(P+B)F = Y$  and since $Y$ is in $\cq_n$ we get 
$\langle BF_j, F_j\rangle = 0 $  for all indices $j,$ Which means that the set of columns in $F$ is not full, and then  the set of columns in $X$ is not full.  
The statement on the dimensions follows from the previous proposition, and the theorem follows. \end{proof}

The theorem has the following immediate corollary.

\begin{corollary} Let $ X$ be in $\cq_n$ and $ X = \sum_{i=1}^k t_iE_i$ be a convex combination of extremal points in $\cq_n$ with sum $X,$ then $k \geq \mathrm{rank}(X)/ \sqrt{n}.$ \end{corollary}

\section{Extremal points in the set of Schur multipliers of norm at most 1}
We have not found any simple characterization of the extremal Schur multipliers of norm at most 1 in $M_{(m,n)}(\bc),$ but we can use the methods from above to get some necessary conditions an extremal Schur multiplier must fulfil. Based on the result in itmem (ii) of Theorem \ref{ScMu} we make the following definition. 
\begin{definition}  Let $X$ be a matrix in $M_{(m,n)}( \bc)$ of rank $r. $  A pair of matrices $L$ in $M_{(r,m)}(\bc)$ with $\|L\|_c = \sqrt{\|S_X\|}$ and $R$ in $M_{(r, n)}(\bc)$ with $\|R\|_c = \sqrt{\|S_X\|}$ is called a Schur factorization of $X$ if $X = L^*R.$ \end{definition}  
It should be remarked, that Schur factorizations are never uniquely determined since we may always replace the pair $(L, R)$ by the pair $(VL, VR)$ for any unitary $V$ in $M_r(\bc).$ And this is not the only way in which different Schur factorizations may occur. 
The item (ii) of Theorem \ref{ScMu} implies that all matrices do have at least one Schur factorization.

\begin{theorem} \label{GenExt} 
Let $X$ in $M_{(m,n)}(\bc)$ be an extremal Schur multiplier of norm 1 and rank $r.$  
\begin{itemize} 
\item[(i)]Let  $X = L^*R$ be a Schur factorization of $X$ then all columns in both $L$ and $R$ are unit vectors.
 \item[(ii)] There exists a factorization $X = FVG$ such that  $F$ is positive of rank $r$ in   $M_m(\bc)$ with  $F^2 $ in $\cq_m,$   $G$ is  positive of rank $r$ in $M_n(\bc) $ with $G^2 $  in $\cq_n$ and $V $ is the matrix of a partial  isometry from  the range of $G$ onto the range of $F.$  
\item[(iii)] The set of columns  $\{ L_1, \dots, L_m, R_1, \dots , R_n \}$ is a full set of vectors in $\bc^r.$
\item[(iv)] The rank  of $X$ is at most $ \sqrt{m+n}.$ 

\end{itemize}
\end{theorem}

\begin{proof}
When $X$ is written in a Schur factorization as $X = L^*R$ and a column, say $R_j$ is not a unit vector in $\bc^r$, then that vector may be written as a  convex combination $R_j = (1-s) \mu + s \nu $ of different vectors $\mu$ and $\nu$ both of norm at most 1. Hence $R$ is a non trivial convex combination $ R = (1-s) R_\mu + s R_\nu$ when the column $R_j$ is replaced by respectively $\mu$ and $\nu.$ 
Since all columns in $L$ span $\bc^r$ we see that both $L^*R_\mu $ and $L^*R_\nu$ are different from $X.$ By item (i) of Theorem  \ref{ScMu} we get that $\|S_{L^*R_\mu}\| \leq 1 $ and $\|S_{L^*R_\nu}\| \leq 1 $ and then we may write $X$ as the non 
trivial convex combination $X = (1-s)L^*R_\mu + sL^*R_\nu,$ which shows that $X$ is not extremal if a column in $R$ is not a unit vector, and item (i) follows.

\smallskip
To see the validity of item (ii) we look at a Schur factorization $X = L^*R$ from where we construct the  polar decompositions $R = UG$ and $L = WF$ with $F$ a positive matrix of rank $r$ such that, by item (i), $F^2 $ is in $\cq_m$ and   $G $ is a positive matrix of rank $r$ with $G^2 $ in $\cq_n.$ Then $U $ is the matrix of  a partial isometry of the range of $G$ onto $\bc^r,$ and $W$ is the matrix of a partial isometry of the range of $F$ onto $\bc^r.$ Then $V:= W^*U $ is the matrix of a patial isometry from the range of $G$ onto the range of $F, $ and item (ii) follows.

\smallskip
Let us suppose that the set of vectors $\{L_1, \dots, L_m, R_1, \dots, R_n\}$ is not full, then there exists a self-adjoint matrix $B$ in $M_r(\bc) $ of norm 1, such that for each index $i$ and each index $j$ we have $$0 = \langle B L_i, L_i\rangle = \langle B R_j, R_j \rangle.$$
As above we have  $I_r \neq I_r + B \geq 0$ and $I_r \neq I_r - B \geq 0. $ Define matrices $A$ and $C$ in $M_{(r,m)}(\bc)$ by $A:= (I_r+B)^{(1/2)} L $ and $C:= (I_r-B)^{(1/2)} L, $ and by the assumption made, we have $\|A\|_c = \|C\|_c  = \|L\|_c = 1.$ Similarly we define matrices $G$ and $H$ in $M_{(r,n)}(\bc)$ by  $G:= (I_r+B)^{(1/2)} R$ and $H:= (I_r-B)^{(1/2)} R, $ and we have $\|G\|_c = \|H\|_c  = \|R\|_c =1.$ Then by item (i) of Theorem  \ref{ScMu} we get that $\|S_{A^*G} \| \leq 1, $ and $\|S_{C^*H}\| \leq 1.$  Since the range of both $L$ and $R$ equals all of $\bc^r$ we have $A^*G = L^*(I_r+B)R \neq X,$ $C^*H = L^*(I_r -B)R \neq X$  and $\frac{1}{2} A^*G +\frac{1}{2} C^*H =X.$ Then   we see that $X$ is not an extremal Schur multiplier of norm 1, if the set of columns is not full, and item (iii) follows.

\smallskip
When $X$ is extremal, the set of all columns is full, so $(m+n ) \geq r^2,$ and the theorem follows. 
\end{proof}
\section{An example} 
Up to now we did not know of a way to obtain extremal points in $\cq_n$ of rank more than 1, but Theorem \ref{ExtQ} can help. That  theorem shows that there are no extremal points in $\cq_3$ except those of rank 1.  In order to get an extremal point in $\cq_4$ of rank 2,  we se from Theorem \ref{ExtQ} combined with the  polar decomposition theorem, that  we  shall  look for a matrix $L$ in $M_{(2,4)}(\bc) $ whose columns  forms a full set of unit vectors in $\bc^2.$ Then an extremal point in $\cq_4$ will be $L^*L.$ As an example we can take $$L := \begin{pmatrix}
1 & 0 & \frac{\sqrt{2}}{2} &  \frac{\sqrt{2}}{2} \\ 
0 & 1 & \frac{\sqrt{2}}{2} &  i\frac{\sqrt{2}}{2} 
\end{pmatrix} \text{ and } L^*L = \begin{pmatrix}
1 & 0 &  \frac{\sqrt{2}}{2} &  \frac{\sqrt{2}}{2} \\
 0 & 1 & \frac{\sqrt{2}}{2} &  i\frac{\sqrt{2}}{2} \\ 
 \frac{\sqrt{2}}{2} &  \frac{\sqrt{2}}{2}  & 1 & \frac{1+i}{2}\\
  \frac{\sqrt{2}}{2} & -i \frac{\sqrt{2}}{2}  &  \frac{1- i}{2} & 1\\  
\end{pmatrix}.$$ Let $k$ be a natural number, then it surprised us, that if we take $k$ unit vectors in $\bc^2$ and supplement $L$  with $k$ extra columns consisting of those unit vectors, we will get  a matrix $\hat L$ in $M_{(2, 4 +k)}(\bc)$ such that  $\hat L^* \hat L$ is extremal in $\cq_{(4 +k)}.$

\end{document}